\documentclass[a4paper,12pt]{amsart}
\usepackage{amsmath,amsthm}
\usepackage{amsfonts,amsmath,amsthm, amssymb, mathrsfs}
\usepackage{float}
\usepackage{euscript,eufrak,verbatim}
\usepackage{graphicx}
\usepackage[usenames]{color}
\usepackage[colorlinks,linkcolor=red,anchorcolor=blue,citecolor=blue]{hyperref}
\usepackage{amsmath}
\usepackage{amsthm}
\usepackage[all]{xy}
\usepackage{graphicx}
\usepackage{amssymb}
\usepackage{enumerate}

\theoremstyle{plain}
\newtheorem{theorem}{Theorem}[section]

\newtheorem{definition}[theorem]{Defintion}
\newtheorem{lemma}[theorem]{Lemma}

\makeatletter 
\@addtoreset{equation}{section}
\numberwithin{equation}{section}

\title{on integral formula  for upper metric mean dimension with potential}

\begin{document}

\title{Entropy Formulae on Feldman-Katok Metric of Random Dynamical Systems}
\author[Yunxiang Xie, Ercai Chen and  Kexiang Yang]{Yunxiang Xie$^1$, Ercai Chen$^1$ and  Kexiang Yang*$^1$}
\address{1.School of Mathematical Sciences and Institute of Mathematics, Nanjing Normal University, Nanjing 210023, Jiangsu, P.R.China}
\email{yxxie20@126.com }
\email{ecchen@njnu.edu.cn}
\email{kxyangs@163.com}

\thanks{}
\subjclass[]{}

\maketitle
\renewcommand{\thefootnote}{}
\footnote{*corresponding author}

\begin{abstract}
	In this paper, we study the Feldman-Katok metric in random dynamical systems and establish corresponding
	fiber topological entropy formula, Brin-Katok local entropy formula and fiber Katok entropy formula by replacing Bowen metric with Feldman-Katok metric. It turns out that the Feldman-Katok metric is also the weakest metric that makes the entropy formulae valid on random dynamical systems.	
\end{abstract}
\maketitle
\section{Introduction}
The setup consists of a probability space ($\Omega$,$\mathcal{F}$,$ \textbf{\textbf{P}} $,$\vartheta$)  together with an invertible ergodic
system $\vartheta$ of a compact metric space $ (X,\mathcal{B}_{X}) $ together with the distance function $d$ and the
Borel $\sigma$-algebra $\mathcal{B}_{X}$.
We assume $\mathcal{F}$ is complete, countably generated, and separatespoints, and so ($\Omega$,$\mathcal{F}$,$ \textbf{\textbf{P}} $) is a Lebesgue space. Let $\mathcal{E}$ =$\Omega$ $\times$ X be  measurable with respect to the product $\sigma$-algebra $\mathcal{F}$ $\times$ $\mathcal{B}$, so the filbers
$\mathcal{E}_{\omega}$=$\{x\in X\colon (\omega,x)\in\mathcal{E}\}$=X, $\forall \omega\in \Omega$.
A continuous bundle random
dynamical system  (abbreviated as $RDS$) $T=(T_{\omega})_{\omega}$ over ($\Omega$,$\mathcal{F}$,$ \textbf{\textbf{P}} $,$\vartheta$) is generated
by mappings $T_{\omega}\colon X\rightarrow X$ with  iterates
$$ T^{n}_{\omega}=
\begin{cases}
id.& \text{$ n=0, $}\\
T_{\vartheta^{n-1}\omega}\circ\cdots T_{\vartheta\omega}\circ\ T_{\omega}& \text{$n \geq 1$,}
\end{cases}$$
where $id$ is an identity mapping, such that  such that $(\omega,x)\mapsto f_{\omega}(x)$ is measurable and the map $x\mapsto f_{\omega}(x)$ is continuous for $\textbf{P}$-$a.e. \omega \in \Omega$.
The map $\Theta $: $\Omega \times X\rightarrow \Omega \times$ $ X $ is defined by
$\Theta (\omega,x)=(\vartheta\omega,T_{\omega}x)$
which is called the skew product transformation.
For $\omega\in\Omega$, $n\in\mathbb{N}$, the $ n$-th  Bowen metric $d^{n}_{\omega}$  on $ X $ is defined by $$d^{n}_{\omega}(x,y)=\max \limits_{0 \leq i \leq n-1}d(T^{i}_{\omega}(x),T^{i}_{\omega}(y)).$$ For   $\forall \varepsilon >0$,  $ B_{n}(\omega,x,\varepsilon) $ denote the open ball with centre $x$ and radius $\varepsilon$ in the metric $d^{n}_{\omega}$, i.e.,
$$ B_{n}(\omega,x,\varepsilon)=\{y\in X\colon d^{n}_{\omega}(x,y)<\varepsilon\}.$$

In \cite{B92}, Bogensch\"{u}tz applied the Bowen metric to  $ RDS $, similarly gave the definition of the corresponding spanning and separated sets and proved the topological entropy of RDS. For a $RDS$, Zhu \cite{Z08,Z09} eatablished the Brin-Katok entropy formula and Katok entropy formula on the Bowen metric.	

By a pair $(X, T)$ we mean a topological dynamical system (abbreviated as $TDS$) where $X$
is a compact metric with the metric $d$ and $T : X\rightarrow X$ is a continuous mapping. In ergodic theory, a fundamental problem is to classify the measure-perserving systems(abbreviated as $MPS$) up to isomorphism. In 1958, Kolmogorov \cite {K58} introduced the concept of  entropy in ergodic theory based on Shannon \cite{S48} entropy, and proved that  entropy is an isomorphic invariant of $MPSs$, then Sinai \cite{S59}  generalized it to the general case.
In ergodic theory, a remarkable result about isomorphic problem is  Ornstein's theorey \cite {O74}, which proves that any two Bernoulli processes of isentropic are isomorphic. The
concept of a finitely determined process plays an important role in his theory, which is based on the hamming distance $\bar{d}_{n}:$
\begin{align*}
\bar d_{n}(x_{0}x_{1}\dots x_{n-1}, y_{0}y_{1}\dots y_{n-1})
=\frac{\lvert\{0\leq i\leq n-1\colon x_{i}\neq y_{i}\}\rvert}{n}.
\end{align*}

In 1943, Kakutani \cite {K43} proposed the equivalent concept between ergodic systems, which is called  Kakutani equivalence. In 1976, by changing the Hamming distance $\bar{d}_{n}$ in Ornsten's theory to the edit distance $\bar{f}_{n}:$
\begin{align*}
\bar f_{n}(x_{0}x_{1}\dots x_{n-1}, y_{0}y_{1}\dots y_{n-1})=1-\frac{k}{n},
\end{align*}
where $ k $ is the largest integer such that there exists
$$0 \leq i_{1}<\dots <i_{k} \leq n-1,0 \leq j_{1}<\dots <j_{k} \leq n-1$$ and $x_{i_{m}}=y_{j_{m}}$ for $m=1,\dots, k,$
Fledman \cite {F76} defined  loose Bernoulli systems, which brings new ideas for the classification of guarantee systems.
In 2017, Kwietniak and Lacka \cite {KM17} introduced the Feldman-Katok (FK) metric as the topological counterpart of edit distance. In 2020, Garc\'{\i}-Ramos and  Kwietniak \cite {G20}  used the FK metric to describe the zero entropy Bernoulli guarantee system. Recently Cai-Li \cite{CL21} established the topological entropy formula of FK metric, Brin-Katok local entropy formula and Katok entropy formula in the case of invariant measures and ergodic measures, and proved that FK metric was the weakest metric to make the topological entropy formula valid.  Nie and Huang \cite{NH22} further studied the restricted-sensitivity of
 FK metirc and mean metric and obtained conditional entropy formulae.

Our purpose in this paper is to extend the results of \cite{CL21} to $ RDS $. In section 2,  we introduce  some basic notions needed in the paper and give entropy formulae for Bowen metric of $RDS$. In section 3, for a continuous bundle $ RDS $, we first prove that the topological entropy defined with  FK metric is equal with Bowen metric.
\begin{theorem}
Let $T$ be a continuous bundle $ RDS $ over $ (\Omega,\mathscr{F},\textbf{\textbf{P}},\vartheta) $. Then for \textbf{P}-$a.e.$ $\omega \in \Omega$
\begin{equation*}
\begin{split}
h_{top}^{(r)}(\omega,X,T)=\overline{h}^{(r)}_{FK}(\omega,X,T)=\underline{h}^{(r)}_{FK}(\omega,X,T).
\end{split}
\end{equation*}
where $ h_{top}^{(r)}(\omega,X,T) $ is the fiber topological entropy of $ T $ with respect to $\omega\in\Omega$.
\end{theorem}

Next, we establish the  Brin-Katok local entropy formula and Katok entropy formula for FK metric.
\begin{theorem}
	Let $T$ be a continuous bundle $ RDS $ over $ (\Omega,\mathscr{F},\textbf{\textbf{P}},\vartheta) $, $\mu\in M^{1}_{\textbf{P}}(\Omega \times X, T)$. Then for $\mu$-$ a.e.$ $(\omega,x)\in \Omega \times X,$
	\begin{align*}
	h^{(r)}_{\mu}(T,\omega,x)&=\lim_{\delta \to 0}\liminf_{n\to \infty}-\frac{1}{n}\log \mu_{\omega}(B_{FK_{n}}(\omega,x,\delta)) \\
	&=\lim_{\delta\to 0}\limsup_{n\to \infty}-\frac{1}{n}\log\mu_{\omega}(B_{FK_{n}}(\omega,x,\delta)).
	\end{align*}
\end{theorem}
\begin{theorem}
	Let $ T $ be a continuous bundle $ RDS $ over  $ (\Omega,\mathscr{F},\textbf{\textbf{P}},\vartheta) $, $\mu\in M^{1}_{\textbf{\textbf{P}}}(\Omega \times X, T)$.  Then  for \textbf{P}-$a.e.$ $\omega \in \Omega$,
	$$h^{(r)}_{\mu}(T)\leq \lim_{\varepsilon \to 0}\liminf_{n\to \infty} \dfrac{1}{n}\log sp_{FK}(\omega,\mu,n,\varepsilon).$$
	If $\mu \in E^{1}_{\textbf{P}}(\Omega \times X, T)$ and $h^{(r)}_{\mu}(T)<\infty $, then
	\begin{align*}
	\begin{split}
	h^{(r)}_{\mu}(T)&= \lim_{\varepsilon \to 0}\liminf_{n\to \infty} \dfrac{1}{n}\log sp_{FK}(\omega,\mu,n,\varepsilon)\\&=\lim_{\varepsilon \to 0}\limsup_{n\to \infty} \dfrac{1}{n}\log sp_{FK}(\omega,\mu,n,\varepsilon),
	\end{split}
	\end{align*}
	where $ sp_{FK}(\omega,\mu,n,\varepsilon) $ denotes  the smallest
	number of any $FK$-$(n, \varepsilon)$-dynamical balls ($ i.e. $ the balls have
	radius $\varepsilon$ in the metric $d^{FK_{n}}_{\omega}$) whose union has $\mu_{\omega}$-measure larger
	than $1-\varepsilon$.
\end{theorem}

\section{Preliminaries}
In this section we will introduce some notions necessary for this paper and  give fiber topological entropy and measure-theoretic entropy formulae with Bowen metric of $RDSs$.
\subsection{Entropies for RDSs}
For RDSs, the definitions of the fiber
topological entropy and fiber measure-theoretic entropies were introduced by \cite{B93,K01,Z08,Z09}.

Let $T$ be a continuous bundle $RDS$ over $ (\Omega,\mathcal{F},\textbf{P},\vartheta) $. For $n \in \mathbb{N},\omega\in \Omega, \varepsilon>0,$ the Feldman-Katok(FK) metric $d^{FK_{n}}_{\omega}$ on $X$ is defined as follows:
for $x,y\in X$, we define an $(\omega,n,\varepsilon)$-match of $x $ and $y$ to be an order preserving (i.e. $\pi(i) < \pi(j)$ whenever $i < j$) bijection $\pi : D(\pi)\rightarrow R(\pi)$
such that $D(\pi), R(\pi)$ are subsets of  $\{0,1,\cdots,n-1 \}$ and for every  $i\in D(\pi)$, it holds  that  $d(T^{i}_{\omega}(x),T^{\pi(i)}_{\omega}(y))<\varepsilon$.
Let
$\lvert \pi \rvert$  denotes the cardinality of $D(\pi)$. Set
\begin{center}
	$\bar{f}_{\omega,n,\varepsilon}(x,y)$=$1-\frac{1}{n}$ max $\{ \lvert \pi \rvert\colon \pi$  is an $ (\omega,n,\delta)$ -match of $ x $ and  $y \}$.
\end{center}	
FK metric of $RDSs$ is
\begin{align}
d_{\omega}^{FK_{n}}(x,y)=\inf \{\varepsilon>0\colon \bar{f}_{\omega,n,\varepsilon}(x,y)<\varepsilon \}.
\end{align}
\begin{definition}
	Let $Z$ be a compact subeset of $X$. A subset $ E\subset Z $ is said to be a $ FK $-($\omega,n,\varepsilon$) spanning set of $ Z $ if $\forall x\in Z, \exists  y\in$ $ E $\, with  $d_{\omega}^{FK_{n}}(x,y)\leq\varepsilon $.
	Let $sp_{FK}(\omega,n,Z,\varepsilon)$ denote the smallest cardinality of any $ FK $-($\omega$,n,$\varepsilon$) spanning set for $ Z $.	
	A set $  F\subset Z $ is said to be a $ FK $-($\omega,n,\varepsilon$) separated set of $Z$ if $\forall  x,y\in F$, $x\neq y$ implies $d_{\omega}^{FK_{n}}(x,y)>\varepsilon $.
	Let $sr_{FK}$($\omega,n,Z,\varepsilon$) denote the largest cardinality of any $ FK $-($\omega,n,\varepsilon$) separated set for $ Z$.	When $Z=X$ we omit the restriction of $Z$.
\end{definition}
It is clear that  $sp_{FK}$($\omega,n,\varepsilon$) $\leq$ $sr_{FK}$($\omega,n,\varepsilon$)  $\leq$ $sp_{FK}$($\omega,n,\frac{\varepsilon}{2}$).
Then the following limits exist$\colon$
\begin{align}
\begin{split}
\overline{h}^{(r)}_{FK}(\omega,X,T)&=\lim_{\varepsilon \to 0}\limsup\limits_{n\to \infty}\frac{1}{n}\log sp_{FK}(\omega,n,\varepsilon)\\&=\lim_{\varepsilon \to 0}\limsup_{n\to \infty}\frac{1}{n}\log sr_{FK}(\omega,n,\varepsilon),\\
\underline{h}^{(r)}_{FK}(\omega,X,T)&=\lim_{\varepsilon \to 0}\liminf\limits_{n\to \infty}\frac{1}{n}\log sp_{FK}(\omega,n,\varepsilon)\\&=\lim_{\varepsilon \to 0}\liminf_{n\to \infty}\frac{1}{n}\log sr_{FK}(\omega,n,\varepsilon).
\end{split}
\end{align}
Replacing $d^{FK_{n}}_{\omega}$ by $d^{n}_{\omega},$ the corresponding notions are $r(\omega,n,Z,\varepsilon)$ and $s(\omega,n,Z,\varepsilon)$.  Then the fiber topological entropy of $T$ with respect to $\omega$ is given by
\begin{equation}\label{Bowen top}
\begin{split}
h_{top}^{(r)}(\omega,X,T)&=\lim_{\varepsilon \to 0}\liminf_{n\to \infty}\frac{1}{n}\log r(\omega,n,\varepsilon)=\lim_{\varepsilon \to 0}\limsup_{n\to \infty}\frac{1}{n}\log r(\omega,n,\varepsilon)  \\
&=\lim_{\varepsilon \to 0}\liminf_{n\to \infty}\frac{1}{n}\log s(\omega,n,\varepsilon)=\lim_{\varepsilon \to 0}\limsup_{n\to \infty}\frac{1}{n}\log s(\omega,n,\varepsilon).
\end{split}
\end{equation}
And the fiber topological entropy of $T$ is given by
$$h^{(r)}_{top}(T)=\int h_{top}^{(r)}(\omega,X,T) d\textbf{P}(\omega).$$

\subsection{Local Entropy of RDSs}
For a continuous bundle $RDS$ $T$, we denote by $\mathcal{P}_{\textbf{P}}(\Omega\times X) $ the space of probability measures on $ \Omega\times X $ having the marginal $ \textbf{P} $ on $\Omega$. We always assume that $\mu$ disintegrates with respect to $ \textbf{P} $, $ i.e. $ there is a family of conditional probabilities $\left \{ \mu_{\omega}  \right \}$ such that
$d\mu(\omega,x)$=$d\mu_{\omega}(x)d\textbf{P}(\omega)$.

A probability measure $\mu$ on $(\Omega \times X,\mathcal{F} \times \mathcal{B})$ is said to be $ T $-invariant if it is invariant under $\Theta$ and has marginal $ \textbf{P} $ on $\Omega$. Furthermore, $\mu$ is said to be $ T $-ergodic if it is ergodic with respect to $\Theta$. Denote by $M_{\textbf{P}}^{1}(\Omega \times X,T)$ the set of all invariant measures of $\Omega \times X$  and  by $E_{\textbf{P}}^{1}(\Omega \times X,T)$  the set of all ergodic measures of $\Omega \times X$. By Bogensch\"{u}tz \cite{B93}, $\mu$
is $\Theta$-invariant if and only if $T_{\omega}\mu_{\omega}=\mu_{\vartheta\omega}$ $\textbf{P}$-$a.e.$.

Let $\mu \in M_{\textbf{P}}^{1}(\Omega \times X,T)$ and $ \zeta  $ be a finite measurable partition of $\Omega \times X $, then the limit
\begin{equation} \label{entropy}
\begin{split}
h^{(r)}_{\mu}(T,\zeta):=\lim\limits_{n \to \infty}\dfrac{1}{n} \int H^{(r)}_{\mu_{\omega}}(\bigvee^{n-1}_{i=0}(T_{\omega}^{i})^{-1}\zeta_{\vartheta^{i}\omega})d\textbf{P}(\omega)
\end{split}
\end{equation}
exists, where $ \zeta_{\vartheta^{i}\omega}$ is the $\vartheta^{i}\omega$-section of $\zeta$ and
$$H^{(r)}_{\mu_{\omega}}(\bigvee\limits^{n-1}_{i=0}(T_{\omega}^{i})^{-1}\zeta_{\vartheta^{i}\omega})=-\sum\limits_{A\in \bigvee\limits_{i=0}^{n-1}(T_{\omega}^{i})^{-1}\xi}\mu_{\omega}(A)\log\mu_{\omega}(A).$$
The number
$$h^{(r)}_{\mu}(T):=\sup\{h^{(r)}_{\mu}(T,\zeta) \mid \zeta\  is \ a \ finite \ measurable \ partition\  of \ \Omega \times X  \} $$
is called the measure-theoretic entropy of $(T,\mu)$.

The classical Brin-Katok entropy formula and Katok entropy formula of the $TDS$ were established by Brin and Katok \cite{BK83,K80}.  Zhu \cite{Z08,Z09} gave a random version of Brin-Katok entropy formula and Katok entropy formula.
Let $T$ be a continuous bundle $RDS$ and $\mu \in M_{\textbf{P}}^{1}(\Omega \times X,T)$, it can be proved that for $\mu$-$a.e.(\omega,x)\in \Omega \times X$, the following equation holds $\colon$
\begin{equation}\label{Bowen Brin}
\begin{split}
h_{\mu}^{(r)}(T,\omega,x)=&\lim_{\varepsilon \to 0}\liminf_{n\to \infty}-\frac{1}{n}\log\mu_{\omega}(B_{n}(\omega,x,\delta))\\
=&\lim_{\varepsilon \to 0}\limsup_{n\to \infty}-\frac{1}{n}\log\mu_{\omega}(B_{n}(\omega,x,\delta)),
\end{split}
\end{equation}
where we call $h_{\mu}^{(r)}(T,\omega,x)$ the local entropy at $(\omega,x)$ with respect to $\mu$. And
\begin{align}
h_{\mu}^{(r)}(T)=\int h_{\mu}^{(r)}(T,\omega,x)d\mu(\omega,x).
\end{align}
Particularly, if $\mu \in E_{\textbf{P}}^{1}(\Omega \times X,T)$,
$ h_{\mu}^{(r)}(T)=h_{\mu}^{(r)}(T,\omega,x)$ for
$\mu$-$a.e.(\omega,x)\in \Omega \times X$.

For $\varepsilon>0,\delta>0,\omega\in \Omega$ and $\mu \in E_{\textbf{P}}^{1}(\Omega \times X,T),$ denote
$$r(\omega,\mu,n,\varepsilon,\delta)=\min\{r(\omega,n,Z,\varepsilon)\colon Z\subset X, \mu_{\omega}(Z)\geq 1-\delta\}.$$
Note that for any fixed $n,\delta$, the map $\varepsilon\mapsto r(\omega,\mu,n,\varepsilon,\delta)$ is monotone decreasing and for any
fixed $n,\varepsilon$, the map $\varepsilon\mapsto r(\omega,\mu,n,\varepsilon,\delta)$ is monotone increasing. Zhu\cite{Z09} proved that for any $\delta \in (0,1)$, if $h_{\mu}^{(r)}(T)< \infty,$ then for $\textbf{P}$-$a.e.$ $\omega \in \Omega$,
\begin{equation}\label{KAO}
\begin{split}
h_{\mu}^{(r)}(T)=&\lim_{\varepsilon \to 0}\limsup_{n\to \infty}\frac{1}{n}\log r(\omega,\mu,n,\varepsilon,\delta)\\
=&\lim_{\varepsilon \to 0}\liminf_{n\to \infty}\frac{1}{n}\log r(\omega,\mu,n,\varepsilon,\delta).
\end{split}
\end{equation}

\section{Entropy formulae for  FK  metric on RDSs}
In this section, we shall extend the results of \cite{CL21} to $ RDSs $.
\subsection{Topological entropy formula}
\begin{lemma}
	Let $T$ be a continuous bundle $RDS$  over $ (\Omega,\mathcal{F},\textbf{P},\vartheta) $, $ n\in \mathbb{N},$ and $x,y\in X $. Then for \textbf{P}-$a.e.$ $\omega \in \Omega$, we have	
	 $d_{\omega}^{FK_{n}}(x,y)\leq d^{n}_{\omega}(x,y)$.
\end{lemma}
\begin{proof}Fix $ x,y \in X $. For any $\varepsilon>d^{n}_{\omega}(x,y)$, then $$\pi\colon \{0,1,\cdots,n-1\}\rightarrow \{0,1,\cdots,n-1\} $$
	is an $(\omega,n,\varepsilon)$-match of $ x $ and $ y $, where $\pi=id.$. Clearly $\bar{f}_{\omega,n,\varepsilon}(x,y)=0$, then
	$$d_{\omega}^{FK_{n}}(x,y)\leq \varepsilon.$$
	Let $\varepsilon \to d^{n}_{\omega}(x,y)$.
\end{proof}	

\begin{theorem}
	Let $T$ be a continuous bundle $ RDS $ over $ (\Omega,\mathcal{F},\textbf{P},\vartheta) $. For $\textbf{P}$-$a.e.$ $\omega \in \Omega$,
	\begin{equation} \label{topological 1}
	\begin{split}
	h_{top}^{(r)}(\omega,X,T)=\overline{h}^{(r)}_{FK}(\omega,X,T)=\underline{h}^{(r)}_{FK}(\omega,X,T).
	\end{split}
	\end{equation}
\end{theorem}
\begin{proof} Let $\mathcal{U}$ be a finite open cover of $ X $ with the Lebesgue number 2$\varepsilon_{0}$. Fix $\omega \in \Omega,n\in \mathbb{N}$. For any 0$<\varepsilon<\varepsilon_{0} $,  Let $ E $ be a $ FK $-($\omega$,n,$\varepsilon$) spanning set with $ \lvert E  \rvert= sp_{FK}(\omega,n,\varepsilon)$. By the
	definitions of $ FK $-($\omega$,n,$\varepsilon$) spanning and $d^{FK_{n}}_{\omega}$, we can get the following relationship$\colon$
	$$X=\bigcup_{x\in E}\bigcup_{k=[(1-\varepsilon)n]}^{n}\mathop{\bigcup_{\pi\colon
			\lvert \pi \rvert=k}}_
	{\pi\ is \ order \ preserving} (T^{i}_{\omega})^{-1} B(T_{\omega}^{\pi(i)}x,\varepsilon).$$
	
	It is obvious that $B(T_{\omega}^{\pi(i)}x,\varepsilon)$ is contained in some element of $\mathcal{U}$. Hence  $\bigcap\limits_{i\in D(\pi)}(T^{i}_{\omega})^{-1}\mathcal{U}$
	is contained in some element of $\bigvee\limits_{i\in D(\pi)}(T^{i}_{\omega})^{-1}\mathcal{U}$.
	We can see that $$\bigvee\limits_{i=0}^{n-1}(T^{i}_{\omega})^{-1}\mathcal{U}=(\bigvee\limits_{i\in D(\pi)}(T^{i}_{\omega})^{-1}\mathcal{U})\bigvee(\bigvee\limits_{i\not\in D(\pi)}(T^{i}_{\omega})^{-1}\mathcal{U}),$$
	and $$\lvert \bigvee\limits_{i\in D(\pi)}(T^{i}_{\omega})^{-1}\mathcal{U} \rvert \leq \lvert \mathcal{U} \rvert ^{n-\pi}.$$
	Hence $\bigcap\limits_{i\in D(\pi)}(T^{i}_{\omega})^{-1}\mathcal{U}$ can  be covered by $ \rvert \mathcal{U} \rvert ^{n-\pi} $ elements of $\bigvee\limits_{i=0}^{n-1}(T^{i}_{\omega})^{-1}\mathcal{U}$.
	Since the number of order preserving bijection $ \pi $ with $\rvert \pi \rvert$= $ k $ is $(C^{k}_{n})^{2}$, it is easy to see that $ X $ can
	be covered by
	\begin{center}
		$\rvert E \rvert \sum\limits_{k=[(1-\varepsilon)n]}^{n}(C^{k}_{n})^{2}\lvert \mathcal{U} \rvert ^{n-k}$
	\end{center}
	elements of $\bigvee\limits_{i=0}^{n-1}(T^{i}_{\omega})^{-1}\mathcal{U}$.
	
	Therefore
	\begin{align*}
	N(\bigvee\limits_{i=0}^{n-1}(T^{i}_{\omega})^{-1}\mathcal{U}) & \leq \rvert E \rvert \sum\limits_{k=[(1-\varepsilon)n]}^{n}(C^{k}_{n})^{2}\lvert \mathcal{U} \rvert ^{n-k} \\
	& \leq sp_{FK}(\omega,n,\varepsilon)\lvert \mathcal{U} \rvert ^{n\varepsilon +1}(n\varepsilon +1)(C_{n}^{[n\varepsilon] +1})^{2}.
	\end{align*}
	Then
	\begin{equation*}
	\begin{aligned}
	\dfrac{1}{n}N(\bigvee\limits_{i=0}^{n-1}(T^{i}_{\omega})^{-1}\mathcal{U})\leq&{\dfrac{\log sp_{FK}(\omega,n,\varepsilon) }{n}+\dfrac{\log (n\varepsilon+1)}{n}}\\
	&+\dfrac{2\log C_{n}^{[n\varepsilon] +1}}{n}+(\varepsilon+\frac{1}{n})\log\lvert \mathcal{U} \rvert.
	\end{aligned}
	\end{equation*}
	By Stirling’s formula,
	$\lim\limits_{n\to \infty}\frac{1}{n}\log C_{n}^{[n\varepsilon] +1}= -(1-\varepsilon)\log(1-\varepsilon)-\varepsilon \log\varepsilon.$
	We have 	
	$$h_{top}^{(r)}(\omega,X,\mathcal{U})\leq\lim_{\varepsilon \to 0}\liminf_{n\to \infty}\frac{1}{n}\log sp_{FK}(\omega,n,\varepsilon).$$
	$ i.e. $
	\begin{align}\label{TOP1}
	h_{top}^{(r)}(\omega,X,T)\leq\underline{h}^{(r)}_{FK}(\omega,X,T).
	\end{align}
	
	On the other hand, since $$d_{\omega}^{FK_{n}}(x,y)\leq d_{\omega}^{n}(x,y),$$
	we have $$sp_{FK}(\omega,n,\varepsilon)\leq r_{n}(\omega,n,\varepsilon),$$
	By (\ref{Bowen top}), it can be obtained that
	\begin{align}\label{TOP2}
	h_{top}^{(r)}(\omega,X,T)\geq\lim_{\varepsilon \to 0}\limsup_{n\to \infty}\frac{1}{n}\log sp_{FK}(\omega,n,\varepsilon)=\overline{h}^{(r)}_{FK}(\omega,X,T).
	\end{align}
	Combining with the fact that $\underline{h}^{(r)}_{FK}(\omega,X,T)\leq \overline{h}^{(r)}_{FK}(\omega,X,T)$, we finish the proof.
\end{proof}

\subsection{Measure-theoretical local entropy formulae of FK metric}
Now we consider the measurable case. First we need some preparation. Recall the definition of edit distance $\bar f_{n}\colon$
\begin{align*}
\bar f_{n}(x_{0}x_{1}\dots x_{n-1}, y_{0}y_{1}\dots y_{n-1})=1-\frac{k}{n},
\end{align*}
where $ k $ is the largest interger such that there exists
$$0 \leq i_{1}<\dots <i_{k} \leq n-1,0 \leq j_{1}<\dots <j_{k} \leq n-1$$ and $x_{i_{m}}=y_{j_{m}}$ for $m=1,\dots, k.$

Let $\mu$ $\in$ $M_{\textbf{P}}^{1}(\Omega \times X,T)$ and $\xi=\left \{A_{1},A_{2} \dots,A_{m}  \right \}$  be a finite partition of $ X $. For any fixed $\omega \in \Omega$, we can identify the elements in $\bigvee\limits_{i=0}^{n-1}(T^{i}_{\omega})^{-1}\xi $ and $\left \{1,2,\dots,m  \right \}^{n}$ by
\begin{align}\label{eaulity}
\bigcap\limits_{i=0}^{n-1}(T^{i}_{\omega})^{-1}A_{t_{i}}=(t_{0},t_{1},\dots,t_{n-1}).
\end{align}
Hence when $ t $$\in$ $\left \{1,2,\dots,m \right \}^{n}$ and $ A, B $$\in$ $\bigvee\limits_{i=0}^{n-1}(T^{i}_{\omega})^{-1}\xi $, we can respectively talk about $\mu_{\omega}$($ t $)
and
$\bar f_{n}(A,B)$.

Next we give the notion of $\pi_{X}$ which is the projection from $\Omega \times X$ onto $ X \colon $
\begin{align*}
\pi_{X}: \  & \Omega \times X \rightarrow X \\
&  (\omega,x)\mapsto x.
\end{align*}
It is obvious that $\pi_{X}$ is measurable. For $ B $$\in$ $\mathcal{B}$, we have
$$\pi_{X} \mu(B)=\mu \circ \pi_{X}^{-1}(B)=\int{\mu_{\omega}(B)}d\textbf{P}(\omega).$$
By the compactness of $ X $, we know that $\pi_{X} \mu$ is a regular measure and
we can construct a finite measurable partition $\eta$ with $\pi_{X} \mu(\partial \eta)$=0, $ i.e. \ $$\mu_{\omega}(\partial \eta)$=0 for \textbf{P}-$a.e.$ $\omega \in \Omega$.

Recall that $B_{n}(\omega,x,\delta)=\{y\in X \colon d^{n}_{\omega}(x,y)< \delta\}  $ is the Bowen ball of $RDS$. Then we replace  $  d_{n}^{\omega} $ by $ d^{FK_{n}}_{\omega}$ and $\bar{f}_{n}$ and denote
$$B_{{FK_{n}}}(\omega,x,\delta)=\{ y\in X \colon d^{FK_{n}}_{\omega}(x,y)< \delta \},$$
$$B_{\bar{f}_{n}}(A,\kappa)=\{ B\in \bigvee\limits_{i=0}^{n-1}(T^{i}_{\omega})^{-1}\xi\colon \bar{f}_{n}(A,B)<\kappa \}.$$

To prove the the Brin-Katok entropy formula for  $d_{\omega}^{FK_{n}}$, Shannon-McMillan-Breiman theorem of $RDS$ needs to be used as a tool, which had been proved in \cite{Z08}.
\begin{lemma}\label{SMB}
	Let $ (X, d) $ be a compact metric space, $T$ a continuous bundle $ RDS $ on $ (X, d)  $ over $ (\Omega,\mathscr{F},\textbf{\textbf{P}},\vartheta) $. For any finite partition $\xi$ of $ X $, if $\mu\in M^{1}_{\textbf{P}}(\Omega \times X, T)$, then we have
	\begin{equation} \label{SMB}
	\begin{split}
	{ \lim_{n \to \infty} -\dfrac{1}{n}\log\mu_{\omega}(\xi^{n}(x)) = h^{(r)}_{\mu}(\xi,\omega,x)}
	\end{split}
	\end{equation}
	where $\xi^{n}(x)$ is the member of the partition $\bigvee\limits_{i=0}^{n-1}(T^{i}_{\omega})^{-1}\xi $ to which $ x $ belongs.
\end{lemma}

Then we can get the following theorem.
\begin{theorem}\label{BK}
	Let $T$ be a continuous bundle $ RDS $ over $ (\Omega,\mathcal{F}, \textbf{P} ,\vartheta) $, $\mu\in M^{1}_{\textbf{P}}(\Omega \times X, T)$, then for $\mu$-$a.e.(\omega,x)\in \Omega \times X,$
	\begin{equation} \label{Brin1}
	\begin{split}
	h^{(r)}_{\mu}(T,\omega,x)&=\lim_{\delta \to 0}\liminf_{n\to \infty}-\frac{1}{n}\log \mu_{\omega}(B_{{FK_{n}}}(\omega,x,\delta)) \\
	&=\lim_{\delta\to 0}\limsup_{n\to \infty}-\frac{1}{n}\log\mu_{\omega}(B_{{FK_{n}}}(\omega,x,\delta)).
	\end{split}
	\end{equation}
\end{theorem}

\begin{proof}	
	We first prove the inequality
	$$h^{(r)}_{\mu}(T,\omega,x) \geq \lim_{\delta \to 0}\limsup_{n\to \infty}-\frac{1}{n}\log \mu_{\omega}(B_{d_{\omega}^{FK_{n}}}(x,\delta)).$$
	Since $d_{\omega}^{FK_{n}}(x,y)\leq d_{\omega}^{n}(x,y)$, by (\ref{Bowen Brin}), we can obtain that
	$$h^{(r)}_{\mu}(T,\omega,x) \geq \lim_{\delta \to 0}\limsup_{n\to \infty}-\frac{1}{n}\log \mu_{\omega}((B_{{FK_{n}}}(\omega,x,\delta)).$$
	
	Now we proceed to prove the  inequality
	$$h^{(r)}_{\mu}(T,\omega,x) \leq \lim_{\delta \to 0}\liminf_{n\to \infty}-\frac{1}{n}\log \mu_{\omega}(B_{{FK_{n}}}(\omega,x,\delta)).$$
	Let $$\mathcal{E}_{0}=\left \{ (\omega,x)\in \Omega \times X \colon h_{\mu}(T,\omega,x)<\infty \right \},$$
	$$\mathcal{E}_{\infty}=\left \{ (\omega,x)\in \Omega \times X \colon h_{\mu}(T,\omega,x)=\infty \right \}.$$
	Without loss of generality, we assume that $\mu(\mathcal{E}_{0})>0, \mu(\mathcal{E}_{\infty})>0$.
	
	Note that $X$ is compact, we can construct a family of increasing finite Borel partitions $\left\{\xi_{i}\right\}_{i=1}^{\infty}$ of $ X $ with $\mu_{\omega}(\partial \xi_{i})=0, \forall i$ and diam$(\xi_{i}) \rightarrow 0, i\rightarrow \infty$.
	
	Given $\varepsilon$ $>$0. Since
	$$ h^{(r)}_{\mu}(\xi,\omega,x) \rightarrow h^{(r)}_{\mu}(T,\omega ,x),a.e.(\omega,x)\in \mathcal{E}_{0},$$
	we can find $\xi$ $\in$ $\left\{\xi_{i}\right\}_{i=1}^{\infty}$ with $\mu$($ A $) $>$ $\mu(\mathcal{E}_{0})-\dfrac{\varepsilon}{2}$,
	where
	\begin{center}
		$ A $=$\left \{ (\omega,x)\in \mathcal{E}_{0} \colon  \lvert   h^{(r)}_{\mu}(\xi,\omega,x) - h^{(r)}_{\mu}(T,\omega, x) \rvert   <  \dfrac{\varepsilon}{2}  \right \}$	.
	\end{center}
	For $\delta >0$, we define
	\begin{equation} \label{eqn2}
	\begin{split}
	U_{\delta}(\xi)=\left \{ (\omega,x)\in \mathcal{E} \colon  B_{{FK_{n}}}(\omega,x,\delta) \backslash \xi(x)\neq \emptyset \right \}.
	\end{split}
	\end{equation}
	It is clear that for $(\omega,x)$ $\in$ $U_{\delta}(\xi)$, the $\delta$-ball about $ x $ on the metric $d^{FK_{n}}_{\omega}$ is not cantained in the element of $\xi$ which $ x $ belongs. On the other hand, we can find that $$U_{\delta}(\xi)=\bigcup\limits_{\omega \in \Omega} (U_{\delta}(\xi))_{\omega} $$  where
	$$(U_{\delta}(\xi))_{\omega}=\left \{x \in X \colon  (\omega,x) \in U_{\delta}(\xi) \right \}, \omega \in \Omega.$$
	
	By the definions above, for \textbf{P}-$a.e.$ $\omega \in \Omega $,  $\bigcap\limits_{\delta > 0}(U_{\delta}(\xi))_{\omega}=\partial \xi.$ Then we can get that
	$$\lim\limits_{\delta \to 0}\mu_{\omega}(U_{\delta}(\xi))_{\omega}=0.$$
	For simplicity, we assume that this convergence is uniform in $\omega$ (otherwise, for $\forall n >0, $ by the Egorov theorem, we can choose $\Omega_{n}\subset \Omega$, with $ \textbf{P} $($\Omega_{n}$) $>$
	1-$\frac{1}{n}$,  such that the above convergence is uniform in $\omega \in \Omega_{n}$). Choose 0 $<$ $\kappa$ $<$ $\varepsilon$ with
	$$ 2\kappa \log\lvert \xi  \rvert -4\kappa \log\kappa-4(1-\kappa)\log(1-\kappa)<\dfrac{\varepsilon}{2}. $$
	Therefore we can find $ \delta \in (0,\dfrac{\kappa}{2})$  such that $\forall  \delta_{0} \leq \delta$, for \textbf{P}-$a.e.$ $\omega \in \Omega,$   $$\mu_{\omega}(U_{\delta_{0}}(\xi))_{\omega} < (\dfrac{\kappa}{4})^{2}$$
	and  $$\mu(U_{\delta_{0}}(\xi))=\int_{\Omega}{(U_{\delta_{0}}(\xi))_{\omega}d\textbf{P}(\omega)} <(\dfrac{\kappa}{4})^{2}.$$
	
	By Birkhoff theorem, $\exists \chi^{*}_{U_{\delta_{0}}(\xi)} \in L^{1}(\Omega \times X,  \mathcal{F} \times \mathcal{B},\mu)$ such that for $\mu$-$a.e.$ $(\omega,x)\in \Omega \times X$,
	\begin{equation} \label{Birkhoff}
	\begin{split}
	\lim_{n \to \infty}\dfrac{1}{n}\sum_{i=0}^{n-1}\chi_{U_{\delta_{0}}(\xi)}(\Theta^{i}(\omega,x))=\chi^{*}_{U_{\delta_{0}}(\xi)}(\omega,x).
	\end{split}
	\end{equation}
	By Chebyshev inequality, we have
	\begin{equation*}
	\begin{aligned}
	\mu(\left \{ (\omega,x)\colon \chi^{*}_{U_{\delta_{0}}(\xi)}(\omega,x) < 1-\dfrac{\kappa}{4} \right \}) & \geq 1-\dfrac{ \int{\chi^{*}_{U_{\delta_{0}}(\xi)}(\omega,x)}d\mu(\omega,x)}{\dfrac{\kappa}{4}} \\
	&= 1-\dfrac{\mu(U_{\delta_{0}}(\xi)) }{\dfrac{\kappa}{4}}\\
	&> 1-\dfrac{\kappa}{4}. \\
	\end{aligned}
	\end{equation*}
	For simplicity as before, we assume that the convegence in (\ref{Birkhoff}) is uniform in $(\omega,x)$. Therefore, when for $ n $ large enough , we can obtain that
	$$\mu(\left \{ (\omega,x)\colon \dfrac{1}{L}\sum_{i=0}^{L-1}\chi_{U_{\delta_{0}}(\xi)}(\Theta^{i}(\omega,x)) < 1-\dfrac{\kappa}{2} ,\forall L\geq n \right \})> 1-\dfrac{\kappa}{2}.$$
	Denote
	 $$E_{L}=\left \{ (\omega,x)\colon \dfrac{1}{L}\sum_{i=0}^{L-1}\chi_{U_{\delta_{0}}(\xi)}(\Theta^{i}(\omega,x)) < 1-\dfrac{\kappa}{2} ,\forall L\geq n \right \}.$$
	Then we can find that $\mu(E_{L})>1-\dfrac{\varepsilon}{2}$.
	Hence $ A $ $\cap$ $E_{L}$ $\subset \mathcal{E}_{0} $ and $\mu(A \cap E_{L})> \mu(\mathcal{E}_{0})-\varepsilon.$ Let
	$$A_{k}=\left \{ (\omega,x)\in A\cap E_{L}\colon k\varepsilon \leq h_{\mu}(T,\omega,x)<(k+1)\varepsilon \right \}.$$
	We have $\bigcup\limits_{k=0}^{\infty}A_{k}=A\cap E_{L}$.
	  Hence there exists some $N_{1}$ such that
	$\mu(\bigcup\limits_{k=0}^{N_{1}}A_{k})>\mu(\mathcal{E}_{0})-\varepsilon$.
	By Lemma \ref{SMB}, we have
	 \begin{align}\label{50}
	\lim_{n \to \infty} -\dfrac{1}{n}\log\mu_{\omega}(\xi^{n}(x)) = h^{(r)}_{\mu}(\xi,\omega,x),  a.e.(\omega,x)\in \Omega\times X .
	\end{align}
	For simplicity , we assume that the convegence in (\ref{50}) is uniform in $(\omega,x)$.
	For $(\dfrac{\varepsilon}{2})^{2}$, $\exists N_{2} \in \mathbb{N},$ then $\forall n>N_{2} $, we have $$-\dfrac{1}{n}\log\mu_{\omega}(\xi^{n}(x))-h^{(r)}_{\mu}(\xi,\omega,x)>(\dfrac{\varepsilon}{2})^{2}.$$
	Then by the Chebyshev inequality, for $ n $ large enough we have
	$$\mu(\left \{ (\omega,x) \in \bigcup_{k=0}^{N_{1}}A_{k} \colon -\dfrac{1}{n}\log\mu_{\omega}(\xi^{n}(x)) > h^{(r)}_{\mu}(\xi,\omega,x)-\dfrac{\varepsilon}{2} \right \})
	> \mu(\bigcup_{k=0}^{N_{1}}A_{K}) -\dfrac{\varepsilon}{2}.$$
	By the  definition of $A_{K}$, we have
	$$\mu(\bigcup_{k=0}^{N_{1}}\left \{ (\omega,x) \in A_{K} \colon -\dfrac{1}{n}\log\mu_{\omega}(\xi^{n}(x)) > h^{(r)}_{\mu}(\xi,\omega,x)-\dfrac{\varepsilon}{2} \right \})
	> \mu(\bigcup_{k=0}^{N_{1}}A_{k}) -\dfrac{\varepsilon}{2}.$$
	Hence we can find $B_{k} \subset A_{k},0\leq k\leq N_{1}$ such that
	$$\mu(\bigcup_{k=0}^{N_{1}}B_{k}) > \mu(\mathcal{E}_{0})-\dfrac{3\varepsilon}{2},$$
	
	Denote the filbers
	$$E_{L,\omega}=\left \{ x\colon \dfrac{1}{L}\sum\limits_{i=0}^{L-1}\chi_{U_{\delta_{0}}(\xi)}(T_{\omega}^{i}(x)) < 1-\dfrac{\kappa}{2} ,\forall L\geq n \right \},$$
	
	$$A_{k,\omega}=\left\{x \colon (\omega,x)\in A_{k} \right\},B_{k,\omega}=\left\{x \colon (\omega,x)\in B_{k} \right\}.$$

		For $ \forall\  n>N_{2},\forall x\in B_{k,\omega}$, we have
		\begin{equation}\label{k}
		\begin{split}
		-\dfrac{1}{n}\log\mu_{\omega}(\xi^{n}(x)) > h^{(r)}_{\mu}(\xi,\omega,x)-\dfrac{\varepsilon}{2}> h^{(r)}_{\mu}(T,\omega,x)-\varepsilon \geq (k-1)\varepsilon.
		\end{split}
		\end{equation}
		The second inequality holds because $x \in A$ and the last inequality holds because $x \in A_{k}.$
		
		Let $N_{0}=\max\{\left\{N_{2},L \right\}\}, n>N_{0}$.

	\textbf{Claim$\colon$} Fix $ \in \Omega.$ For $x \in B_{k,\omega}$, we have
	$$ B_{{FK_{n}}} (\omega,x,\delta_{0}) \subset \bigcup\limits_{t \colon \bar{f}_{n}(t,\xi^{n}(x))<\kappa}t .$$

	\textbf{Proof of claim$\colon$} Let $ y\in B_{{FK_{n}}} (\omega,x,\delta_{0}) .$ We only need to prove that there exists $t=(t_{0},\cdots,t_{n-1})$ such that $y\in\bigcap\limits_{i=0}^{n-1} (T_{\omega}^{i})^{-1}A_{t_{i}}$
	and $\bar{f}_{n}(t,\xi^{n}(x))<\kappa.$
	
	By the definition of $d^{FK_{n}}_{\omega}$, for fixed $\omega \in \Omega, n \in \mathbb{N},x,y \in X ,$ there exists an $ (\omega,n,\delta_{0})$-match of $ x $ and $ y $ with
	$$ \lvert \pi  \rvert =\lvert D(\pi)  \rvert >n(1-\delta_{0}).$$
	Since  $x \in E_{L,\omega}$, we have
	$$ \lvert \left \{ 0\leq i\leq n-1\colon T_{\omega}^{i}(x)\notin (U_{\delta_{0}}(\xi))_{\omega} \right \}  \rvert>n(1-\dfrac{\kappa}{2}).$$
	Note that $$  \lvert  D(\pi)  \rvert >n(1-\delta_{0}) >n(1-\dfrac{\kappa}{2}), $$ we have
	$$ \lvert \left \{ i\in  D(\pi) \colon T_{\omega}^{i}(x)\notin (U_{\delta_{0}}(\xi))_{\omega} \right \}  \rvert>n(1-\dfrac{\kappa}{2}).$$
	Denote $$D_{1}=\left \{ i\in  D(\pi) \colon T_{\omega}^{i}(x)\notin (U_{\delta_{0}}(\xi))_{\omega} \right \} .$$
	Clearly
$D_{1} \subset D(\pi)$ and for $\forall j $$\in D_{1},$   $$T_{\omega}^{j}(x)\notin (U_{\delta_{0}}(\xi))_{\omega},$$  by  (\ref{eqn2})  we can obtain that $$  B(T_{\omega}^{j},\delta_{0})\subset \xi(T_{\omega}^{j}x). $$
	For $\forall j\in D_{1},$   $d(T_{\omega}^{j}(x),T_{\omega}^{\pi(j)}(y))<\delta_{0}$, we have
	$$T_{\omega}^{\pi(j)}(y)\in B(T_{\omega}^{j},\delta_{0}).$$ Hence $T_{\omega}^{\pi(j)}(y)\in \xi(T_{\omega}^{j}x).$ It follows that
	$$\bar{f}_{n}(\xi^{n}(y),\xi^{n}(x))\leq 1-\dfrac{ \lvert  D_{1} \rvert}{n}<\kappa.$$
	The claim is proved.
	
	Next for $ \omega \in \Omega,$ we estimate the $\mu_{\omega}$-measure of the set
	$$\{x\in B_{k,\omega} \colon \mu_{\omega}(B_{{FK_{n}}} (\omega,x,\delta_{0}))>e^{-n(k-2)\varepsilon} \}.$$
	Note that the number of elements in $B_{\bar{f}_{n}}(\xi^{n}(x),\kappa)$ is not more than
	$(C_{n}^{[n\kappa]})^{2}|\xi|^{n\kappa }$, then we can see that
	\begin{equation*}
	\begin{aligned}
	& \{x\in B_{k,\omega} \colon \mu_{\omega}(B_{{FK_{n}}} (\omega,x,\delta_{0}))>e^{-n(k-2)\varepsilon} \} \\
	\subset \ & \{x\in B_{k,\omega} \colon \mu_{\omega}(\bigcup\limits_{t \colon \bar{f}_{n}(t,\xi^{n}(x))<\kappa}t)>e^{-n(k-2)\varepsilon} \}  \\
	\subset \ & \{x\in B_{k,\omega} \colon \exists \ t \in B_{\bar{f}_{n}}(\xi^{n}(x),\kappa), s.t. \mu_{\omega}(t)>\frac{e^{-n(k-2)\varepsilon}}{(C_{n}^{[n\kappa]})^{2}|\xi|^{n\kappa }} \} \\
	\subset \ &  \bigcup_{t\bigcap  B_{k,\omega}\neq \emptyset } \{t \colon \exists\  t^{'} \in B_{\bar{f}_{n}}(t,\kappa), s.t. \mu_{\omega}(t^{'})>\frac{e^{-n(k-2)\varepsilon}}{(C_{n}^{[n\kappa]})^{2}|\xi|^{n\kappa }} \}
	\end{aligned}
	\end{equation*}
	For the last set, the number of such $t$ is not more than $$ e^{n(k-2)\varepsilon}((C_{n}^{[n\kappa]})^{2}|\xi|^{n\kappa })^{2}. $$
	We only need to estimate the $\mu_{\omega}$ of a single $t$. Let $x\in t \cap B_{k,\omega}$. Since $t \in \bigvee\limits_{i=1}^{n-1}(T_{\omega}^{i})^{-1}\xi,$
	it follows that
	$ t \in \xi^{n}(x),$ hence $t=\xi^{n}(x)$. By (\ref{k})  we have
	$$-\dfrac{1}{n}\log\mu_{\omega}(\xi^{n}(x)) > (k-1)\varepsilon.$$
	That is $\mu_{\omega}(t)<e^{-n(k-1)\varepsilon}$. Hence
	\begin{equation*}
	\begin{aligned}
	& \mu_{\omega}\{ x\in B_{k,\omega} \colon \mu_{\omega}(B_{{FK_{n}}} (\omega,x,\delta_{0}))>e^{-n(k-2)\varepsilon}  \}\\
	<& \ e^{-n(k-1)\varepsilon} e^{n(k-2)\varepsilon}((C_{n}^{[n\kappa]+1})^{2}|\xi|^{n\kappa })^{2} \\
	=&\ e^{-n\varepsilon}((C_{n}^{[n\kappa]})^{2}|\xi|^{n\kappa })^{2}.
	\end{aligned}
	\end{equation*}
	By Stirling’s formula, $$\lim_{n\to \infty}\frac{1}{n}\log C_{n}^{[n\varepsilon] }= -(1-\varepsilon)\log(1-\varepsilon)-\varepsilon \log\varepsilon. $$
	Clearly
	$$\lim_{n\to \infty} (C_{n}^{[n\kappa]})^{4}|\xi|^{2n\kappa +2}=e^{2n\kappa \log\lvert \xi  \rvert -4n\kappa \log \kappa-4(1-\kappa) \log(1-\kappa)}<e^{\frac{n\varepsilon}{2}}. $$
	Then  $\sum\limits_{n} e^{-n\varepsilon}(C_{n}^{[n\kappa]})^{4}|\xi|^{2n\kappa +2}$ is convergent. By the Borel-Cantelli lemma, we have
	$$\liminf_{n\to \infty}-\dfrac{\log\mu_{\omega}(B_{d^{FK_{n}}_{\omega}}(x,\delta_{0}))}{n}\geq (k-2)\varepsilon>h^{(r)}_{\mu}(T,\omega,x)-3 \varepsilon, a.e. (\omega,x)\in B_{k}.$$
	Hence we have
	$$\liminf_{n\to \infty}-\dfrac{\log\mu_{\omega}(B_{{FK_{n}}} (\omega,x,\delta_{0}))}{n}>h^{(r)}_{\mu}(T,\omega,x)-3\varepsilon, a.e. (\omega,x)\in \bigcup^{N}_{k=0}B_{k}.$$
	Note that $\mu(\bigcup\limits^{N}_{k=0}B_{k})> \mu(\mathcal{E}_{0})-\dfrac{3\varepsilon}{2}$, then by the arbitrariness of $\varepsilon$, we can get that
	$$ \lim_{\delta \to 0}\liminf_{n\to \infty}-\frac{1}{n}\log \mu_{\omega}(B_{FK_{n}}(\omega,x,\delta))\geq h^{(r)}_{\mu}(T,\omega,x),a.e.(\omega,x)\in \mathcal{E}_{0}.$$
	Similarly we can prove
	$$ \lim_{\delta \to 0}\liminf_{n\to \infty}-\frac{1}{n}\log \mu_{\omega}(B_{FK_{n}}(\omega,x,\delta))\geq h^{(r)}_{\mu}(T,\omega,x) ,a.e.(\omega,x)\in \mathcal{E}_{\infty}.$$
	Therefore $$ \lim_{\delta \to 0}\liminf_{n\to \infty}-\frac{1}{n}\log \mu_{\omega}(B_{FK_{n}}(\omega,x,\delta))\geq h^{(r)}_{\mu}(T,\omega,x) ,a.e.(\omega,x)\in \Omega \times X.$$
	The proof of Theorem $\ref{BK}$ has been finished.
\end{proof}
Next we similarly introduce the measure-theoretic version of $sp_{FK}(\omega,n,\varepsilon)$.
\begin{lemma}
	If $0\leq p_{i}\leq 1,\forall\  0\leq i\leq 1$, and $\sum\limits_{i=1}^{n}p_{i}=1$, then
	$$\sum\limits_{i=1}^{n}(-p_{i}\log p_{i})\leq \log n.$$
\end{lemma}
\begin{definition}
	Let $ T $ be a continuous bundle $ RD $S over ($\Omega$,$\mathcal{F}$,$ \textbf{P} $,$\vartheta$), $\mu\in M^{1}_{\textbf{P}}(\Omega \times X, T)$.
	Denote
	\begin{align*}
	sp_{FK}(\omega,\mu,n,\varepsilon)=\min\{sp_{FK}(\omega,n,Z,\varepsilon)\colon Z\subset X,\mu_{\omega}(Z)>1-\varepsilon \}.
	\end{align*}
\end{definition}	

\begin{theorem}
	Let $ T $ be a continuous bundle $ RDS $ over $ (\Omega,\mathcal{F}, \textbf{P} ,\vartheta) $, $\mu\in M^{1}_{\textbf{P}}(\Omega \times X, T)$.  Then  for \textbf{P}-$a.e.$ $\omega \in \Omega$, we have 	
	$$h^{(r)}_{\mu}(T)\leq \lim_{\varepsilon \to 0}\liminf_{n\to \infty} \dfrac{1}{n}\log sp_{FK}(\omega,\mu,n,\varepsilon).$$
	If $\mu \in E^{1}_{\textbf{P}}(\Omega \times X, T)$ and $h^{(r)}_{\mu}(T)< \infty$, then
	\begin{equation} \label{Katok 1}
	\begin{split}
	h^{(r)}_{\mu}(T)=&\lim_{\varepsilon \to 0}\liminf_{n\to \infty} \dfrac{1}{n}\log sp_{FK}(\omega,\mu,n,\varepsilon)\\
	=&\lim_{\varepsilon \to 0}\limsup_{n\to \infty} \dfrac{1}{n}\log sp_{FK}(\omega,\mu,n,\varepsilon).
	\end{split}
	\end{equation}
\end{theorem}

\begin{proof}
	Note that $d^{FK_{n}}_{\omega}(x,y)\leq d^n_{\omega}(x,y)$ and when $\mu \in E^{1}_{\textbf{P}}(\Omega \times X, T)$, by (\ref{KAO}), for \textbf{P}-$a.e.$ $\omega \in \Omega$, we have
	\begin{equation} \label{Katok 2}
	\begin{split}
	h^{(r)}_{\mu}(T)\geq\lim_{\varepsilon \to 0}\limsup_{n\to \infty} \dfrac{1}{n}\log sp_{FK}(\omega,\mu,n,\varepsilon),
	\end{split}
	\end{equation}
	
	Then we only need to prove that
	when  $\mu\in M^{1}_{\textbf{P}}(\Omega \times X, T)$, we have
	$$h^{(r)}_{\mu}(T)\leq \lim_{\varepsilon \to 0}\liminf_{n\to \infty} \dfrac{1}{n}\log sp_{FK}(\omega,\mu,n,\varepsilon).$$
	Given a finite Borel partition $ \eta $ of $ X $ and $\delta>0$. Let $  \eta=\{ A_{1},\cdots, A_{k}\} $. Note that $\{\Omega \times A_{1},\cdots, \Omega \times A_{k}\} $ is a special finite partition of $ \Omega \times X $, thus  we need to show
	$$h^{(r)}_{\mu}(T,\eta)\leq \lim_{\varepsilon \to 0}\liminf_{n\to \infty} \dfrac{1}{n}\log sp_{FK}(\omega,\mu,n,\varepsilon)+2\delta.$$
	Take $0<\kappa<\frac{1}{2}$ with
	$$-4\kappa\log2\kappa-2(1-2\kappa)\log(1-2\kappa)+(\kappa^{2}+3\kappa)\log(k+1)<\delta.$$
	
	Note that $\pi_{X}\mu$ is normal, then for each $ A_{i} $ we can find a close set $ B_{i} $ with $\pi_{X}\mu(A_{i} \setminus B_{i} ) < \frac{\kappa^{2}}{k}$. Hence  we can construct a  finite Borel partition $ \xi= \{ B_{1},\cdots, B_{k}, B_{k+1}\}$ of $ X $ such that  $ B_{i} $ is closed for
	$1\leq  i\leq k,$ $ \pi_{X}\mu( B_{k+1})<\kappa^{2} ,$ and
	$$ h^{(r)}_{\mu}(T,\eta)<h^{(r)}_{\mu}(T,\xi)+ \delta. $$
	
	Let $ B $=$\bigcup\limits_{i=1}^k B_{i}, b=\min\limits_{1\leq i<j\leq k}d(B_{i},B_{j}).$ Denote $$ B_{\Omega}=\Omega \times B= \{\Omega \times B_{1},\cdots, \Omega \times B_{k}   \}.$$ It is easily to see that $b>0$ and $\mu(B_{\Omega})> 1-\kappa^{2}$.
	
	Let $0< \varepsilon < \frac{b}{2}, n\in \mathbb{N}.$ $\forall \omega \in \Omega$, by the definition of $ sp_{FK}(\omega,\mu,n,\varepsilon),$  there exists $Z\in X $ such that  $ \mu_{\omega}(\bigcup\limits_{i=1}^{m(n)} B_{FK_{n}}(\omega,x_{i},\varepsilon)) >1-\varepsilon$, where $m(n)= sp_{FK}(\omega,\mu,n,\varepsilon).$ For simplicity, we denote
	$$F_{n}=\{(\omega,x)\colon  \omega\in \Omega,     x\in \bigcup\limits_{i=1}^{m(n)}B_{FK_{n}}(\omega,x_{i},\varepsilon) \} $$
	and $$E_{n}=\{(\omega,x)\colon \dfrac{1}{n}\sum\limits_{i=0}^{n-1}\chi_{B}(\Theta^{i}(\omega,x))\leq 1-\kappa \}.$$
	We can get that $\mu(E_{n})< \kappa$. Put $ W_{n}=B_{\Omega}\cap F_{n}\cap E_{n}^{c} $ and
	denote the filber
	$$W_{n,\omega}=\{x\colon  (\omega,x)\in W_{n}\}.$$
	then we can get that $\mu_{\omega}(W_{n,\omega})>1-\kappa^{2}-
	\kappa-\varepsilon$.
	For $z\in W_{n,\omega}$, we have
	$$ \dfrac{1}{n}\sum\limits_{i=0}^{n-1}\chi_{B}(T_{\omega}^{i}z)> 1-\kappa .$$
	
	\textbf{Claim:} For any fixed $ \omega \in \Omega, 1\leq i\leq k$,
	$$\lvert \{ A\in \bigvee_{j=0}^{n-1}(T_{\omega}^{i})^{-1}\xi \colon A\cap B_{FK_{n}}(\omega,x_{i},\varepsilon)\cap W_{n,\omega} \neq \emptyset    \}  \rvert\leq(C_{n}^{[n(2\kappa+2\varepsilon)]})^{2}\lvert   \xi\rvert^{n(2\kappa+2\varepsilon)}.$$
	
	\textbf{Proof of claim.}
	Let
	$ A_{1},A_{2}\in \bigvee\limits_{j=0}^{n-1}(T_{\omega}^{i})^{-1}\xi$  with
	$$ A_{1}\cap B_{FK_{n}}(\omega,x_{i},\varepsilon)\cap W_{n,\omega} \neq \emptyset, A_{2}\cap B_{d^{FK_{n}}_{\omega}}(x_{i},\varepsilon)\cap W_{n,\omega} \neq \emptyset.$$
	For $\forall x\in A_{1}, y\in A_{2} $, we have $d^{FK_{n}}_{\omega}(x,y) < 2\varepsilon$.
	Then  there exists an $ (\omega,n,2\varepsilon) $-match $\pi$ of $x$ and $y$ with $$\lvert D(\pi) \rvert > n(1-2\varepsilon) .$$
	\\
	Denote
	$$D_{x}=\{0\leq j\leq n-1 \colon T_{\omega}^{j}x\in B \},D_{y}=\{0\leq j\leq n-1 \colon T_{\omega}^{j}y\in B \}.$$
	Let $$D^{'}=\pi^{-1}(\pi(D(\pi)\cap D_{x})\cap D_{y}).$$ For simplicity, we assume that $ D^{'}$ is not empty. It is easy to see that $ D^{'}\subset D(\pi)$ and $\lvert  D^{'} \rvert > n(1-2\kappa-2\varepsilon)$.
	For every $j \in D^{'}, $ $$d(T_{\omega}^{j}x,T_{\omega}^{\pi(j)}y)< 2\varepsilon<b $$
	and $T_{\omega}^{j}x,T_{\omega}^{\pi(j)}y \in B$. Hence $T_{\omega}^{j}x,T_{\omega}^{\pi(j)}y$ must lie one of the the same set in $\{B_{1},\cdots,B_{k}\}$.
	It follows that
	$$\bar{f}_{n}(A_{1},A_{2})<2\kappa +2\varepsilon.$$
	Note that the number of $ A $ satisfying
	$$\bar{f}_{n}(A_{1},A)<2\kappa +2\varepsilon$$
	is not more than
	$$ (C_{n}^{[n(2\kappa+2\varepsilon)]})^{2}\lvert   \xi\rvert^{n(2\kappa+2\varepsilon)}.$$
	The proof of the claim is finished.
	
	Now we can estimate $h^{(r)}_{\mu}(T,\xi)$. First we  estimate $H^{(r)}_{\mu_{\omega}}(\bigvee\limits^{n-1}_{i=0}(T_{\omega}^{i})^{-1}\xi)$.
	\begin{equation*}
	\begin{aligned}
	&H^{(r)}_{\mu_{\omega}}(\bigvee\limits^{n-1}_{i=0}(T_{\omega}^{i})^{-1}\xi)\\
	\leq&  H^{(r)}_{\mu_{\omega}}(\bigvee\limits^{n-1}_{i=0}(T_{\omega}^{i})^{-1}\xi\vee\{W_{n,\omega},X\setminus W_{n,\omega} \})\\
	\leq& \mu_{\omega}( W_{n,\omega})\log( \lvert  \{ A\colon A\in \bigvee\limits^{n-1}_{i=0}(T_{\omega}^{i})^{-1}\xi,A\cap W_{n,\omega} \neq \emptyset\}  \rvert)-\mu_{\omega}( W_{n,\omega}) \log \mu_{\omega}( W_{n,\omega})\\
	+&\mu_{\omega}( W^{c}_{n,\omega})\log( \lvert  \{ A\colon A\in \bigvee\limits^{n-1}_{i=0}(T_{\omega}^{i})^{-1}\xi,A\cap W_{n,\omega} \neq \emptyset\}  \rvert)-\mu_{\omega}( W^{c}_{n,\omega})\log \mu_{\omega}( W^{c}_{n,\omega}) \\
	\leq &\mu_{\omega}( W_{n,\omega})\log(\sum_{i=1}^{m(n)}  \lvert  \{ A\colon A\in \bigvee\limits^{n-1}_{i=0}(T_{\omega}^{i})^{-1}\xi,A\cap B_{FK_{n}}(\omega,x_{i},\varepsilon)\cap W_{n,\omega} \neq \emptyset\}  \rvert)\\
	 +&(\kappa^{2}+\kappa+\varepsilon)n\log (k+1)+\log2.
	\end{aligned}
	\end{equation*}
	Thus
	\begin{equation*}
	\begin{aligned}
	h^{(r)}_{\mu}(T,\xi)=&\lim\limits_{n \to \infty}\dfrac{1}{n} \int H^{(r)}_{\mu_{\omega}}(\bigvee^{n-1}_{i=0}(T_{\omega}^{i})^{-1}\xi)d\textbf{P}(\omega)\\
	\leq&\lim_{\varepsilon \to 0}\liminf_{n\to \infty}\dfrac{1}{n}\log sp_{FK}(\omega\,\mu,n,\varepsilon)+2\limsup_{n\to \infty}\dfrac{1}{n}\log C_{n}^{[n(2\kappa+2\varepsilon)]}\\
	&+(\kappa^{2}+3\kappa+3\varepsilon)\log(k+1).
	\end{aligned}
	\end{equation*}
	By Stirling’s formula, let $\varepsilon \to 0$ , we can obtain
	\begin{equation*}
	\begin{aligned}
	h^{(r)}_{\mu}(T,\xi)\leq&\lim_{\varepsilon \to 0}\liminf_{n\to \infty}\dfrac{1}{n}\log sp_{FK} (\omega,\mu,n,\varepsilon)\\
	&-4\kappa\log2\kappa-2(1-2\kappa)\log(1-2\kappa)+(\kappa^{2}+3\kappa)\log(k+1)\\
	<&\lim_{\varepsilon \to 0}\liminf_{n\to \infty}\dfrac{1}{n}\log sp_{FK} (\omega,\mu,n,\varepsilon)+\delta.
	\end{aligned}
	\end{equation*}
	Hence
	$$h^{(r)}_{\mu}(T,\eta)<h^{(r)}_{\mu}(T,\xi)+\delta<\lim_{\varepsilon \to 0}\liminf_{n\to \infty}\dfrac{1}{n}\log sp_{FK}(\omega,\mu,n,\varepsilon)+2\delta.$$
	By the arbitrariness of $\delta$, we have
	\begin{align}h^{(r)}_{\mu}(T)\leq \lim_{\varepsilon \to 0}\liminf_{n\to \infty}\dfrac{1}{n}\log sp_{FK}(\omega,\mu,n,\varepsilon).
	\end{align}	
\end{proof}

\section*{Acknowledgements}

The work was supported by NNSF of China
(11671208 and 11431012). We would like to express our gratitude to Tianyuan Mathematical Center in Southwest China, Sichuan University and Southwest Jiaotong University for
their support and hospitality. The third author is also supported by Postgraduate Research \& Practice Innovation Program of Jiangsu Province (No. KYCX22\_1529).

\end{document}